\documentclass {elsart}
\usepackage [applemac] {inputenc} 
\usepackage [frenchb] {babel}
\usepackage{amssymb}
\usepackage{amsmath}
\usepackage{enumerate}

\newcommand{\K}{\mathbb{K}}
\newcommand{\N}{\mathbb{N}}

\newcommand{\J}{\mathbb{J}}
\newtheorem{definition}{{\bf Definition}}[section]
\newtheorem{theorem}[definition]{{\bf Theorem}}

\newtheorem{corollary}[definition]{{\bf Corollary}}
\newtheorem{proposition}[definition]{\noindent {\bf Proposition}}

\newtheorem{remark}[definition]{\noindent {\bf Remark}}

\newtheorem{lemma}[definition]{\noindent {\bf Lemma}}
\def\endproof{\hfill {\kern 6pt\penalty 500
\raise -0pt\hbox{\vrule \vbox to5pt {\hrule width 5pt
\vfill\hrule}\vrule}}}
\def\centerpicture #1 by #2 (#3){\leavevmode
        \vbox to #2{
        \hrule width #1 height 0pt depth 0pt
        \vfill
        \special{pictfile #3}}}
\begin{document}
\begin{frontmatter}
\title{Incidence structures and Stone-Priestley duality}
\author{ Mohamed Bekkali\corauthref{cor1}}
\corauth[cor1]{Corresponding author.}
\address{ D\'epartement de Math\'ematiques, Facult\'e des Sciences et Techniques,  Universit\'e Sidi Mohamed Ben Abdellah,  B. P. 2202 Sa\"{\i}ss-F\`es, Maroc}

\ead{bekka@menara.ma}

\author{Maurice Pouzet\corauthref{cor2}}
\corauth[cor2]{Supported by INTAS}
\address {PCS, Universit\'e Claude-Bernard Lyon1,
 Domaine de Gerland -b\^at. Recherche [B], 50 avenue Tony-Garnier, F$69365$ Lyon cedex 07, France}
\ead{pouzet@univ-lyon1.fr }

\author {Driss Zhani}
\address{D\'epartement de Math\'ematiques, Facult\'e des Sciences et Techniques,  Universit\'e Sidi Mohamed Ben Abdellah, B. P. 2202 Sa\"{\i}ss-F\`es, Maroc}
\ead{drisszhani@hotmail}
\date{\today}
\begin{center}
Mailbox
\end{center}
\begin{abstract} We observe that if $R:=(I,\rho, J)$ is an incidence 
structure, viewed as a matrix,  then the topological closure of the set of columns is the Stone space of 
the  Boolean algebra  generated by the rows. As a consequence, we 
obtain that  the topological closure of the collection of principal 
initial segments of  a poset $P$ is the Stone space of the Boolean algebra 
$Tailalg (P)$
generated by the collection of principal final segments of $P$, the 
so-called {\it tail-algebra of $P$}. Similar   
results concerning Priestley spaces and distributive 
lattices are  given. A generalization to incidence structures valued by abstract algebras is considered. 
\end{abstract}

\end{frontmatter}

\section{Introduction}
Basic objects  of data analysis  are  matrices  filled with $0$-$1$ coefficients.  The analysis of their structural properties   single out   an  important class of  properties:  those for which the orders on the rows and  on  the columns are irrelevant, that is properties which, once true for a given $m\times n$-matrix $A$, are also  true for all matrices obtained from $A$ by arbitrary permutations of its rows and  columns.  Each row $A_{i,-}$ of $A$ determines a subset of the index set of colums, and dually for the columns.  Simple manipulations of rows and columns translate to set theoretic operations on the corresponding subsets.  For  example, the set made of intersections of columns $A_{-,j}$ ($j:= 1,\dots, n$) is a lattice,  every  matrix obtained from $A$ by permuting its rows and columns  yields  an isomorphic lattice, and moreover, this  lattice is dually isomorphic to the lattice made of intersections of rows. In terms of incidence relations, this   later statement  says that the dual of the Galois lattice of an incidence relation is isomorphic to the Galois lattice of the dual incidence relation, a result at the very heart  of formal concept analysis.  An other result of the  same flavor,  is the following:
\begin{theorem}The number of distinct Boolean combinations made with
the rows of a finite matrix  consisting of $0$ and $1$ is
$2^c$ where $c$ is the number of distinct columns of the matrix.  
\end {theorem}
The proof is immediate.  For, let $A$ be such  a $m\times n$ matrix. Set $I:= \{1, \dots, m\}$, $J:=\{1, \dots, n\}$, $\mathcal R:= \{ k(A_{i,-}) : i\in I \}$, where $k(A_{i, -}):= \{j\in J: A(i,j)=1\}$ and  let $ \mathcal B$ be the set  of Boolean combinations made of members of $\mathcal R$. Say that two indices
$j',j''\in J$ are {\it equivalent} if  the corresponding columns are
identical. There are $c$ equivalence classes. Moreover, if $U$ is an equivalence
class then clearly for every $j'\in 
J\setminus U$ and $j''\in U$ there is some $i\in I$ such that 
$A(i,j')\not = A(i,j')$. It follows then that  $U$ is the intersection of 
 subsets $X$ such that $U\subseteq  X \subseteq J$ and either $X$ or 
$J\setminus X$ belongs to $\mathcal{R}$. Thus, $U\in \mathcal B$ and, hence, $\mathcal {B}$ consists of all unions of equivalence classes. The result 
follows.\\

This result, as basic as it may  seem, might be  as old as the  notion of Boolean algebra. In this paper,  we   first  point out that the extension of this result  to infinite matrices,   captures the essence of Stone  and Priestley dualities alike.   Next, we give an illustration   with the notion of Tail algebra. Finally,  we consider  an extension to matrices with coefficients into an abstract algebra.

For convenience,  we present our results  in terms of incidence structures, rather than matrices.

\subsection{Stone-Priestley duality for incidence structures}

 An {\it incidence structure} is a triple $R:=(I, \rho, J)$ 
where $\rho$ is a relation from a set $I$ to  a set $J$, identified with a subset of the 
cartesian product $I\times J$. For 
$(i,j)$ in $I\times J$, set $R^{-1}(j):=\{i\in I: (i,j)\in \rho \}$ 
and $R(i):= \{j\in J: (i,j)\in \rho \}$; also set  
$\mathcal{C}_R:= \{ R^{-1}(j): j\in J\}$ and  $\mathcal {R}_{R}:=\{R(i): 
i\in I\}$. 

 \noindent Throughout this paper,  
$J$  will denote  a  non-empty set and $\mathfrak{P}(J)$ its power set.  Viewing $\mathfrak{P}(J)$ as a bounded lattice (resp. a Boolean algebra), we denote by $\mathcal
{L}(R)$ (resp.  $\mathcal{B}(R)$), 
the bounded sublattice (resp. the Boolean subalgebra) of $\mathfrak{P}(J)$ 
generated by $\mathcal {R}_{R}$. 
 
\noindent Hence, $\mathcal{L}(R)$ is the smallest collection of 
subsets of $J$ such that (1) $\emptyset$, $J$  and every member of 
$\mathcal{R}_{R}$ belongs to $\mathcal{L}(R)$; (2) $I'\cup I''$ and  $I'\cap I''$ belong to 
$\mathcal{L}(R)$ whenever $I'$ and $I''$ belong to $\mathcal{L}(R)$. Similarly, $\mathcal{B}(R)$ is 
the smallest collection of subsets of $J$ such that (1) every member of 
$\mathcal{L}_{R}$ belongs to $\mathcal{B}(R)$; (2) $I'\cup I''$ and  $I\setminus I'$ belong to 
$\mathcal{B}(R)$ whenever $I'$ and $I''$ belong to $\mathcal{B}(R)$. 
 
\noindent  Identifying  
 $\mathfrak{P}(I)$ with $2^{I}$, we may view it  as a topological space.  A basis of open sets consists of  subsets  of  the form 
$O(F,G):=\{X\in \mathfrak{P}(I): F\subseteq X $ and $ G\cap X=\emptyset \}$, where $F, G$ are 
finite subsets of $I$. Let $\overline {\mathcal{C}}_R$ denotes the topological closure of $\mathcal{C}_R$ in
$\mathfrak{P}(I)$. 
 Recall that a  compact  totally disconnected space is called a {\it Stone space}, whereas a {\it Priestley space} is a set together with a  topology and an ordering which is compact and totally order disconnected. We will use only the fact that closed subspaces of $\mathfrak{P}(I)$, with the inclusion order possibly added, are of this form\cite{prie}.
 
{\bf Example}

\begin{enumerate}
\item  Let  $R:=(E, \in, 
\mathfrak{P}(E))$ where $E$ is a set. Then 
$\mathcal {C}_{R}=\mathfrak{P}(E)$ (thus, it is closed), whereas $\mathcal{R}_R:=\{\{X: x\in X\in 
\mathfrak{P}(E)\} : x\in E\}$.  One may show that  $\mathcal{B}(R)$ is the free Boolean algebra 
generated by $E$ and $\mathcal{C}_R$ is its Stone space.
\item  Replacing  $R$ by its 
dual $R^{-1}:= (\mathfrak{P}(E),\ni ,E)$,  $\mathcal{B}(R)$ is then the power 
set $\mathfrak{P}(E)$ and $\overline {\mathcal{C}_{{R}}}$ is the set 
$\beta(E)$ of ultrafilters on $E$, the Cech-Stone compactification of $E$.
\end{enumerate}
The content of these examples holds in a more general setting.
\begin{theorem} \label{duality}
 The set $\overline {\mathcal{C}}_R$, endowed with the topology induced by  
the powerset $\mathfrak{P}(I)$, is homeomorphic to the Stone space of 
$\mathcal{B}(R)$. With the order of inclusion added, $\overline 
{\mathcal{C}}_R$ is isomorphic to the 
Priestley space of $\mathcal{L}(R)$.
\end{theorem}

\begin{proof} Let $\varphi:\mathfrak{P}(\mathfrak{P}(J)) \rightarrow 
\mathfrak{P}(I)$ be  defined by  $\varphi(\mathcal {B}):= \{i\in 
I : R(i)\in \mathcal {B}\}$ for all $\mathcal{B}\subseteq 
\mathfrak{P}(J)$. Looking at $\mathfrak{P}(\mathfrak{P}(J))$ and  $\mathfrak{P}(I)$ 
 as topological spaces,  we can see that  $\varphi$ is continuous, whereas viewing these sets
 as Boolean  algebras, $\varphi$  is a Boolean homomorphism and, in particular, it 
preserves the ordering.  Let $\mathcal Spec(L)$ be the collection of prime  filters 
of $L$, the {\it spectrum} of $L$, for $L\in \{\mathcal{L}(R),\mathcal{B}(R)\}$.  
We claim that  $\varphi$ induces an isomorphism 
from $Spec(L)$ onto $\overline {\mathcal{C}}_R$ (this isomorphism being a topological one if $L=\mathcal{B}(R)$, a
topological and ordered one if 
$L=\mathcal{L}(R)$).  To this end, let $e: J\rightarrow 
\mathfrak{P}(L)$ be the map defined by setting $e(j):= \{X \in L: j\in 
X\}$. Clearly, $e(j)\in \mathcal Spec(L)$ 
and $\varphi (e(j))= R^{-1}(j)$ for every $j\in J$. It follows 
that $\mathcal{C}_R\subseteq \varphi' (\mathcal Spec(L))$ (where $\varphi' 
(\mathcal Spec(L)):=\{\varphi(\mathcal U): \mathcal U\in \mathcal 
Spec(L)\}$); hence 
$\overline {
\mathcal{C}}_R\subseteq 
\varphi'(\mathcal Spec(L))$, since $\varphi$ is continuous.   Next, 
$
\varphi' (\mathcal Spec(L)) \subseteq \overline {\mathcal{C}}_R$. Indeed, let 
$\mathcal U\in \mathcal Spec(L)$ and let $\mathcal {O}$ be an open set in $\mathfrak{P}(I)$ 
containing $\varphi(\mathcal U)$; without loss of generality, we may 
assume that  $\mathcal {O}=O(F,G):=\{X\in \mathcal 
{P}(I): F\subseteq X $ and $ G\cap X=\emptyset \}$, where $F, G$ are 
finite subsets of $I$. Then, since $\mathcal {U}$ is a prime filter, 
the set $H:= \cap \{R(i): i\in F\}\setminus \cup \{R(i): i\in G\}$ 
belongs to $\mathcal {U}$ and  hence $H$  is non-empty. For  $j\in H$, 
$R^{-1}(j)\in \mathcal {O}$; this proves that $\varphi(\mathcal U)\in  
\overline {\mathcal{C}}_R$. To finish up  the proof, note that $\varphi$ 
is $1-1$ on $ \mathcal Spec(L)$. 
\end {proof}

Every Boolean algebra is of the form $\mathcal B(R)$.  Indeed,  if B is a Boleean algebra,  set $R:= (B, \in , S)$ where  $S$ is the Stone space made of ultrafilters of $B$.  In the next section, we introduce a proper class of Boolean algebras. 
Later on, in  Section \ref{three},  we shall delineate the exact content of 
Theorem \ref{duality} in terms of abstract algebras.

\section{Tail algebras and Tail lattices}
Let $P$ be a poset.  For $x\in P$, the {\it principal final}
(resp. {\it initial}) {\it segment generated by} $x$ is 
$\uparrow x:=\{y\in P : x\leq 
y\}$  (resp. $\downarrow x:=\{y\in P : y\leq x \}$). Set $up(P):=\{\uparrow x: x\in P \}$ and 
$down(P):=\{\downarrow x: x\in P\}$.

The {\it tail  algebra of $P$} is the 
subalgebra $Tailalg( P)$ of the Boolean algebra $(\mathfrak{P}(P), \cap, 
\cup , \setminus , \emptyset , P)$ generated by 
$up(P)$. According to J.D.Monk (\cite{monk}Chap. 
2, p.40), this 
notion is due to G.Brenner. Denote by $Taillat (P)$ the bounded sublattice of 
$(\mathfrak{P}(P), \cup , \cap , \emptyset , P)$, generated by 
$up(P)$. 
Taking $R:=(P,\leq, P)$ in  Theorem \ref {duality} we have: 
\begin {theorem}\label{PPS}

The topological closure $\overline {down(P)}$ of $down(P)$ in $\mathfrak{P}(P)$ is 
homeomorphic to the Stone space of $Tailalg (P)$. With the order of 
inclusion added, $\overline{down(P)}$ it is isomorphic to the Priestley space of 
$Taillat (P)$.
\end{theorem}

The topological closure of $down(P)$ points to interesting 
collections of subsets of $P$. 

A subset $I$ of $P$ is an {\it initial 
segment} (or is {\it closed downward}) if $x\leq 
y$ and $y\in I$ imply $x\in I$; if in addition $I$ is non-empty and up-directed (that is every pair 
$x, y\in I$ has an upper bound $z\in I$), this is an {\it ideal}. For example, each principal initial
segment is an ideal. Let $X$ be a subset of $P$.  We set $\downarrow X:=\{y\in P: y\leq x$ for
some
$x\in X\}$; this set is the least initial segment containing $X$, we say that it is {\it generated}  by $X$;  if $X$ contains only one element $x$,  we will continue to denote it  by $\downarrow x$ instead of $\downarrow \{x\}$. We set $X^{-}:=\bigcap\{\downarrow x : x\in X\}$, the 
 set of lower bounds of $X$. An initial segment $I$ of $P$ is {\it finitely
generated} if $I=\downarrow X$ for some finite subset $X$ of $P$. We denote respectively by
$\bold {I}(P)$, $\bold{I}_{<\omega}(P)$ and  
$\mathcal {J}(P)$  
the collection of initial segments, finitely generated initial segments, and ideals of $P$  ordered by inclusion.
The {\it  dual} de $P$ is the poset obtained from $P$ by reversing the order; we denote it by $P^{*}$
(instead of  $P^{-1}$). Using the above definitions, a subset
 which is respectively an initial segment, a finitely generated initial segment or an ideal of  $P^*$ will
be called a {\it final segment}, a {\it  finitely generated final segment} or a {\it filter}
of
$P$. We denote  by
$\bold {F}(P)$, $\bold{F}_{<\omega}(P)$, and 
$\mathcal {F}(P)$ respectively, the collections of these sets, ordered by inclusion.  Also, denote respectively by $\uparrow X$ and  
$X^{+}$  the least final segment containing
$X$ and the set of upper bounds of $X$ in $P$.

The topological condition of closedness translates in an order theoretic one as shown by the following lemma.
  \begin{lemma}\label{closed}
 Let $X$ be a subset of $P$.   Then $X\in\overline{down(P)}$ if and only if $F^{+}\setminus \uparrow
G\not= \emptyset$ for every finite subsets $F\subseteq X, G\subseteq P\setminus X$.
\end{lemma}
\begin{proof} Clearly, for every $F, G\in \mathfrak {P}(P)$, we  have 
$O(F,G) \cap down(P) \not = \emptyset$ if and only if $F^{+}\setminus \uparrow
G\not= \emptyset$. Let $X\in \mathfrak{P}(P)$. Since the  $O(F,G)$ 's, 
for all finite subsets $F\subseteq X, G\subseteq P\setminus X$, form a basis of neighborhoods of $X$, the
lemma follows.\end{proof}

As an  immediate corollary,  we have:
\begin{corollary} $\emptyset \not \in \overline {down(P)}\Longleftrightarrow P\in
\bold{F}_{<\omega}(P)$
\end{corollary}

We also have:
\begin{corollary}\label{ideals1}
$down(P)\subseteq \mathcal {J}(P)\subseteq \overline {down (P)}\setminus\{\emptyset\}$. In
particular,   the topological closures in $\mathfrak{P}(P)$ of
$down(P)$  and $\mathcal {J}(P)$ are the same. 
\end{corollary}
\begin{proof}
Trivially $down(P) \subseteq \mathcal {J}(P)$; thus, it suffices to check that 
$\mathcal {J}(P)\subseteq \overline {down (P)}$.  Let $I\in \mathcal {J}(P)$. Let  $F, G$ be two finite subsets
of 
$ I$ and $ P\setminus I$ respectively. Since $F$ is finite  and
$I$ is an ideal, there is some $x\in I$ such that $F\subseteq 
\downarrow x$, and since $I$ is an initial segment, $\uparrow G\subseteq P\setminus I$. Hence $x\in
F^{+}\setminus \uparrow G$. According to Lemma \ref  {closed}, $I\in \overline{down (P)}$. This finishes the proof of Corollary \ref{ideals1}.

\end{proof}
\begin{remark}  If $P\not \in \bf {F}_{<\omega}(P)$, $\overline {down(P)}$ is isomorphic (as a Priestley space) to  $\overline {down(P')}$, where $P'$ is the poset obtained from $P$ by adding a least element.
\end{remark}

A poset $P$ is {\it up-closed} if every  intersection of two 
members of $up(P)$ is a finite union (possibly empty) of members of $up(P)$.

\begin{proposition}\label{ideals}
The following properties for a poset $P$ are equivalent:
\begin{enumerate}[{(a)}] 
\item  $\mathcal {J}(P)\cup 
\{\emptyset \}$ is closed for the product topology;
\item  $\mathcal {J}(P)=\overline 
{down(P)}\setminus \{\emptyset\}$;
\item $P$ is up-closed;
\item $\bold {F}_{<\omega}(P)$ is a meet-semilattice; 
\item $Taillat(P)=\bold{F}_{<\omega}(P)\cup \{P\}$.
\end {enumerate}  
\end{proposition}
\begin{proof}
$(a)\Rightarrow (c)$ Suppose that $\mathcal {J}(P)\cup 
\{\emptyset \}$ is closed. Let $x,y \in P$. We prove that $\uparrow x\cap\uparrow y$ is finitely
generated. We may suppose
$\uparrow x\cap\uparrow y\not =
\emptyset$ (otherwise the condition is satisfied).   The set $O(\{x,y\}, \emptyset)\cap \mathcal {J}(P)$
 is included  into $\cup\{O( \{z\}, \emptyset ) : z\in \uparrow x \cap \uparrow y\}$. Being 
compact, from our  hypothesis, it is included into a finite union $\bigcup\{O( \{z\}, \emptyset ) : 
z\in G\}$. Hence $\uparrow x\cap 
\uparrow  y= \cup \{\uparrow z : z\in G \}$ as required.

$(c)\Rightarrow (e)$ Trivial.

$(e) \Rightarrow (d)$ Trivial.

$(d) \Rightarrow (c)$ Let $x,y\in P$. Assuming  that $\bold {F}_{<\omega}(P)$ is a meet-semilattice, it contains an element 
$Z:= \uparrow x \wedge \uparrow y$.   This element is equal to $Z':= \uparrow x \cap \uparrow y$.  
Indeed,  we have obviously $Z\subseteq  Z'$.   If the inclusion was strict, then $Z\cup \uparrow z$, 
for $z\in Z'\setminus Z$, would be  an element $\bold {F}_{<\omega}(P)$ in  between, which is impossible.
 Hence $Z'\in \bold {F}_{<\omega}(P)$, proving $(c)$. 

$(c)\Rightarrow (b)$ Suppose that $P$ is up-closed.  Let $I\in 
\overline {down(P)}\setminus \{\emptyset\}$. Since $\bold I(P)$ is closed in $\mathfrak{P}(P)$,
$I\in 
\bold I(P)$. Let $x,y\in I$;  according to our hypothesis  
 $\uparrow x \cap \uparrow y= \cup \{\uparrow z: z\in G \}$ for some finite 
subset $G$ of $P$ . Necessarily  $G\cap 
I\not =\emptyset$, otherwise $I\in O(\{x, y\}, G)$, whereas $O(\{x, y\}, 
G)\cap \downarrow (P)=\emptyset$, contradicting $I\in 
\overline {down(P)}$. Every $z\in  O(\{x, y\}, G)\cap 
I$ is an upper bound of $x,y$ in $I$ , proving $I\in \mathcal {J}(P)$. 

$(b)\Rightarrow (a)$ Trivial.
    \end{proof}\\
\begin{corollary} \label{ideauxclos}
The following properties  for a poset $P$ are equivalent:
\begin{enumerate}
 \item $\mathcal J(P)$ is closed in $\mathfrak{P}(P)$;
\item $P\in \bold{F}_{<\omega}(P)$ and $P$ is up-closed.
\end{enumerate}
\end{corollary}
As an immediate consequence:

\begin{fact} If $P$ is a join-semilattice with a least element then $\mathcal
J(P)$ is closed for the product topology.
\end{fact}
Let us recall  that if $L$ is a  join-semilattice, an element $x\in L$ is {\it join-irreducible} (resp.
 {\it join-prime}) if it is distinct from the least  element $0$, if any, and if 
$x=a\vee b$ implies
$x=a$ or $x=b$ (resp. $x\leq a\vee b$ implies
$x\leq a$ or $x\leq b$) see \cite{gratzer} where  $0$ is allowed. We denote $\J_{irr}(L)$ ( 
resp. $\J_{pri}(L)$)  the set of join-irreducible (resp. join-prime) members of $P$.
We recall that $\J_{pri}(L)\subseteq \J_{irr}(L)$; the equality holds provided that $L$
 is a distributive lattice. It also holds if $L=\bold {I}_{<\omega}(P)$. Indeed:

\begin{fact}For  an arbitrary poset $P$, we have:
\begin{equation}\label{eqirr}
 {\J}_{irr}(\bold I_{<\omega} (P))=\J_{pri}(\bold I_{<\omega} (P))= down (P)
 \end{equation}

\begin{equation}\label{eqid}
\J_{irr}(\bold{I}( P))= \mathcal {J}(P) 
\end{equation}
\end{fact}
\begin{fact} For a poset $L$, the following properties are equivalent.
\begin{itemize}
\item $L$ is isomorphic to $\bold {I}_{<\omega}(P)$ for some poset $P$; 

\item $L$ is a join-semilattice with a least element  in which  every element is a finite  join of primes. \end{itemize}\end{fact}

We have the following characterization.
\begin{proposition}\label{taillat2}
Let $L$ be a bounded distributive lattice. The following properties are equivalent:

\begin{enumerate}[(a)]

\item Every element  of $L$ is a finite join  of join-irreducible elements; 
\item $L$ is isomorphic to $\bold F_{<\omega}(P)$ for some poset $P$ with $P$ up-closed and $P\in \bold F_{<\omega}(P)$;  
 \item  The Priestley space $Spec(L)$ is isomorphic to  $\mathcal {J}(P)$, for some  poset $P$ with $P$ up-closed and $P\in \bold F_{<\omega}(P)$.
 \end {enumerate}  
\end{proposition} 
\begin{proof}

$(a)\Rightarrow (b)$ Set $Q:={\J}_{irr}(L)$ and let $\varphi :L\rightarrow \mathfrak {P}(Q)$ be defined by  $\varphi(x):= \{y\in Q: y\leq x\}$.  It  is well known that $\varphi$  is an homomorphism of bounded lattices, provided that $L$ is itself a bounded distributive lattice.   Moreover, $\varphi$ is an isomorphism from  $L$ onto $\bold I_{<\omega}(Q)$ if  and only if $L$ satisfies hypothesis  $(a)$. Hence, under this condition, $L$ is isomorphic to $\bold I_{<\omega}(Q)$, that is to $\bold F_{<\omega}(P)$, where $P:= Q^{*}$. To conclude, it suffices  to check that 
$P$ is up-closed and $P\in \bold F_{<\omega}(P)$. As for the proof of  $(d)\Rightarrow  (c)$ in Proposition \ref{ideals},  this simply follows from the fact that $\bold F_{<\omega}(P)$ is a meet-semilattice, with a top element. 

$(b)\Rightarrow (c)$ Assuming that $(b)$ holds, $L$ is isomorphic to $Taillat(P)$. From Theorem \ref{PPS},  $Spec (L)$ is isomorphic to $\overline {down (P)}$, which according to Corollary  \ref  {ideauxclos},  is isomorphic  to $\mathcal J(P)$.

$(c)\Rightarrow (a)$ This implication follows in the same way that $(b) \Rightarrow (c)$ from  Theorem \ref{PPS} and  Corollary  \ref  {ideauxclos}. 

\end{proof}

{\bf Examples} \begin{enumerate}
\item An  exemple from the theory of relations illustrates  Theorem \ref{PPS} and Corollary \ref{ideauxclos}.

First, recall that a {\it relational structure} is a pair $R:=(E, (\rho_i)_i\in I)$, where
for each $i\in I$,  $\rho_i$ is a $n_i$-{\it ary relation} on $E$ (that is a subset of  $E^{n_i}$) and
$n_i$ is a non-negative integer; the family $\mu:=(n_i)_i\in I$ is called the {\it signature} of $R$. An {\it induced substructure} of $R$ is any relational structure of the form $R_{\restriction {F}}:=(F,
(\rho_i\cap F^{n_i})_{i\in I})$. One may define  the notion of  {\it relational isomorphism} and then
the notion of {\it embeddability} between relational structures with the same signature (e.g. $R$
embeds into $R'$ if $R$ is isomorphic to an induced substructure of $R'$). According to Fra{\"\i}ss\'e \cite{fraissetr}, the
{\it age } 
of a relational structure $R$ is  the collection $\mathcal{A}(R)$ of its finite 
induced substructures, considered up to isomorphism. A  first order sentence (in the language associated with the signature $\mu$) is {\it universal} whenever it is equivalent to  a
sentence of the form
$\forall x_1\cdots
\forall x_n \varphi(x_1, \dots, x_n) $ where
$\varphi(x_1, \dots, x_n)$ is a formula built with the variable $x_1,\dots, x_n$, the logical
connectives $\neg$, $\vee$, $\wedge$ and predicates $=$,  $\rho_i, i\in I$.

 Let $\Omega_{\mu}$ be the set of finite  relational structures with signature $\mu$, these
structures being  considered up to  isomorphism and ordered by embeddability. If $\mu$ is finite, then  $\Omega_{\mu}$ is a ranked poset with a least element which is up-closed, hence from Corollary  \ref{ideauxclos},  the set ${\mathcal J}(\Omega_{\mu})$ is a closed subset of $\mathfrak((\Omega_{\mu})$. The use  of ${\mathcal J}(\Omega_{\mu})$ is justified by the following:  
 \begin{proposition}
 ${\mathcal J}(\Omega_{\mu})$ is the set 
of  ages of relational structures with signature $\mu$ and its 
   dual  is the Boolean algebra made of Boolean combinations of  universal sentences,
considered up to elementary equivalence.
    \end{proposition}
 
 Thus from Theorem \ref{PPS} the  tail algebra
$Tailalg(\Omega_{\mu})$ provides an alternative  description of the algebra made of Boolean combinations of  universal sentences ( see \cite {pouzetcc} 
\cite  {pouzetsobranisa} for more details on this correspondence).

If  $\mu$ is constant and equal  to $1$ ($I:= \{0, \dots, k-1\}$ and
$n_i=1$ for $i\in I$) then   $\Omega_{\mu}$ can be identified to the direct product of 
$2^k$ copies of the chain $\omega$ of non-negative integers. From this it follows that  ${\mathcal J}(\Omega_{\mu})$ is homeomorphic to the ordinal $\omega^{2^{k}}+1$, equipped with the interval topology. If $\mu$ contains some integer larger than $1$, then $\Omega_{\mu}$ embeds as a subposet the set $[\omega]^{<\omega}$ of finite subsets of $\omega$ hence ${\mathcal J}(\Omega_{\mu})$ embeds the Cantor space $\mathfrak{P}(\omega)$.

\item  
Algebraic lattices provide an other variety of examples. Indeed, if  $P$ is a join-semilattice with  $0$, then $\mathcal {J}(P)$ is an algebraic lattice, in fact the lattice of closed sets  
of an  algebraic closure system, and every algebraic lattice has this 
form \cite {gratzer}. The fact that algebraic lattices are Priestley spaces leads to  
interesting results, see eg  \cite {mislove},  and questions. Typical examples of 
algebraic lattices are  $\bold {I}(P)$ and  $\bold {F}(P)$, ordered by 
inclusion. For an example, the compact elements of $\bold{F}(P)$ are the finitely 
generated final segments of $P$, from which follows that  $\bold {F}(P) $ is isomorphic to 
$\mathcal{J}(\bold{F}_{<\omega}(P))$. The importance of ${\bold F}(P)$ steems from  the following notion.

The  {\it free Boolean algebra generated  by  a poset $P$} is  a Boolean  algebra containing  $P$  in such a way that every order-preserving map from $P$ into a Boolean algebra $B$ extends to a homomorphism from this Boolean algebra into $B$; up-to an isomorphism fixing $P$ pointwise, this Boolean algebra is unique, it is denoted by $FB(P)$. An explicit description may be given in terms of  tail algebra.

\begin{proposition}
The free-Boolean algebra generated by $P$ is isomorphic to the  tail algebra $Tailalg (\bold{F}_{ 
<\omega}(P))$. \end{proposition}

 For this, note that the
incidence structures $(P, \in, 
\bold{F}_{<{\omega}}(P)$) and $(\bold{F}_{<{\omega}}(P), \subseteq, 
\bold{F}_{<{\omega}}(P))$ yield the same Boolean algebra.  
 
Tail algebras generated by
join-semilattices  with  $0$ are studied in  \cite {bekkali}, \cite {zhani}.
\end{enumerate}
\section {A generalization to valued incidence structures}\label{three}

Let $E$ be a set.   For $n\in \N^*:= \N\setminus
\{0\}$, a map $f: E^n
\rightarrow E$ is an
$n$-{\it ary operation} on $E$, whereas a subset 
$\rho \subseteq E^{n}$ is an $n$-{ary relation} on $E$.
Denote by
$\mathcal O^{(n)}$ (resp.$\mathcal R^{n}$) the set of
$n$-ary operations (resp. relations)  on $E$ and set 
$\mathcal O:=\bigcup \{\mathcal O^{(n)}: n\in N^* \}$ (resp 
$\mathcal R:= \bigcup \{\mathcal R^{(n)}:
n\in N^* \}$). For  $n, i\in N^*$  with $i\leq n$, define the   $i^{th}$
$n$-{\it ary projection}
$e^{n}_{i}$ by setting $e^{n}_{i}(x_{1},\dots,x_{n}):= x_{i}$ for all  $x_{1},\dots,x_{n}\in E$ and set 
$\mathcal P:= \{e^{n}_{i}: i,n \in \N^*\}$. An operation $f\in  \mathcal O$ is {\it constant} if it takes a single value, it 
 is\ {\it idempotent} provided
$f(x,\dots,x)= x$ for all $x\in E$. We denote by  $\mathcal C$ (resp. $\mathcal I$) the set of  constant, (resp. idempotent)  operations on $E$. 

Let $m,n\in \N^*$,
$f \in \mathcal O^{(m)}$ and $\rho \in \mathcal R^{(n)}$.    Then $f$ {\it preserves}
$\rho$
 if: 
\begin{equation}
\small {(x_{1,1}, \dots, x_{1,n})\in
\rho, \dots,
(x_{m,1}, \dots, x_{m,n})\in \rho \Longrightarrow (f(x_{1,1}, \dots, x_{m,1}),\dots,f(x_{1,n}, \dots,
x_{m,n}))\in \rho}
\end{equation}
for every $m\times n$ matrix
$X:= (x_{i,j})_{i=1,\ldots,m \atop {j=1, \ldots ,n}}$ of elements of $E$.  Let $g\in \mathcal O^{(n)}$, 
then $f$ {\it commutes} with $g$ if $f$ preserves the $n+1$-ary relation $$\rho_g:=
 \{(x_1, \dots, x_n, g(x_1, \dots, x_n)): (x_1, \dots, x_n)\in E^{n}\}.$$
 A
{\it universal algebra} (resp. a {\it relational structure}) on
$E$ is a pair
$(E,
\mathcal F)$ where $\mathcal F$ is a subset of $\mathcal O$ (resp. of $\mathcal R$\footnote{Contrarily to the definition given in the previous section, we do not require that $\mathcal F$ is a family of members of $\mathcal R$}.
 Powers of such
structures, subalgebras, and  homomorphisms are easy to define. For  example, if $\K:= (E,
\mathcal F)$ is an algebra, then a subset $\rho \subseteq E^{n}$ induces a  subalgebra of $\K^{n}$ if
every $g\in \mathcal F$ preserves $\rho$; also 
$f:\K^{m}\rightarrow
\K$ is an homomorphism if
$f$ commutes with  
$\rho_{g}$ for all $g\in \mathcal F$; equivalently, $\rho_f$ is a subalgebra of $\K^{m+1}$. 
 
 Let $\K:= (E, \mathcal F)$ be a universal algebra, $I, J$ be two sets  and $A$ be a map from the
direct product
$I\times J$ into 
$\K$.  Given $i\in I$, let  
$A_{i,-}:J\rightarrow \K$ be defined by $A_{i,-}(j):=A(i,j)$; 
similarly, given $j\in J$ let $A_{-,j}:I\rightarrow \K$ be defined by 
$A_{-,j}(i):=A(i,j)$. Set $\mathcal{R}_{A}:=\{A_{i,-} :i\in I\}$ 
and $\mathcal{C}_{A}:=\{A_{-,j} :j\in J\}$. 
 Looking  at $A$ as a matrix, $\mathcal{R}_{A}$ and $\mathcal{C}_{A}$ 
 are the sets of rows and columns of $A$. Let $B$ be the subalgebra of $\K^J$ generated by 
$\mathcal{R}_{A}$, let $Hom(B, \K)$ be the set of homomorphisms from $B$ 
into $\K$     and let $\overline{\mathcal{C}_{A}}$ be the topological 
closure of $\mathcal{C}_{A}$ in $\K^I$ equipped with the product 
topology, $\K$ being equipped with the discrete topology.

 We say that $\K$ is {\it projectively trivial} if every homomorphism $f$ from 
every 
subalgebra $L$ of a finite power $\K^n$ into $\K$ is induced by a 
projection, that is $f(x_{1},\dots, x_{n})=x_{j}$ for some $j$ and 
all $(x_{1},\dots, x_{n})\in L$. \\
Theorem  \ref {duality}  is a consequence of the following statement.
\begin{proposition} 
\begin{enumerate}[(a)] 

\item If  $\K$ is a finite projectively trivial  algebra, then  
 $Hom(B, \K)$ equipped with the topology induced by the product 
 topology on $\K^B$ is homeomorphic to the closure $\overline{\mathcal{C}_{A}}$ of ${\mathcal{C}_{A}}$.
\item The $2$-element 
Boolean algebra and the $2$-element bounded lattice are projectively trivial . 
\end{enumerate}
\end{proposition}
\begin{proof}
(a) We start with $\K$ arbitrary. Let $Y$ be such that 
$\mathcal{R}_{A}\subseteq Y \subseteq \K^J$.
Let $\varphi: \K^{Y}\rightarrow \K^I$, setting 
$\varphi(h)(i):=h(A_{i, -})$ for all $h\in \K^{Y}$, $i\in I$. Let $e: 
J\rightarrow \K^{Y}$ be defined by  
$e(j)(h):=h(j)$ for all $j\in J$, $h\in \K^J$.\\
{\bf Claim 1} $\varphi(e(j))=A_{-,j}$\\
Let $D:=\{g\in \K^I$ such that  $g(i)=g(i')$ whenever 
$A_{i,-}=A_{i',-}\}$. Let $Z$ 
be a subset 
of $\K^Y$ containing the image of $I$ by $e$ and let $\varphi'(Z)$ be the 
image of $Z$ by $\varphi$. \\
{\bf Claim 2} $\mathcal{C}_{A}\subseteq \varphi'(Z) \subseteq D$ \\
Next, we set $Y:=B$ and $Z:=Hom(B,\K)$ .\\
{\bf Claim 3}  $Im e\subseteq Hom(B, \K)$ and $\varphi$ is $1-1$ 
on $Hom(B, \K)$.\\
\begin{proof}Let $h_{1},h_{2}\in Hom(B,\K)$ so that 
$\varphi(h_{1})=\varphi(h_{2})$. This says  
$h_{1}(A_{i,-})=h_{2}(A_{i,-})$ for all $i\in I$, meaning that 
$h_{1},h_{2}$ coincide on $\mathcal{R}_{A}$. Being morphisms,  
$h_{1},h_{2}$ coincide on the algebra generated by $\mathcal{R}_{A}$, 
this algebra is $B$  
hence 
 $h_{1}=h_{2}$ proving that $\varphi$ is $1-1$. 
\end{proof}
  
Let $Im\varphi_{B}$ be the image of $Hom(B,\K)$ under $\varphi$.  
  
  {\bf Claim 4} If $\K$ is finite then $\overline {\mathcal{C}_{A}}\subseteq Im \varphi_{B}$.
 
 \begin{proof}
First, $Hom(B, \K)$ is closed in $\K^B$;  next $\varphi $ is 
continuous. Since $\K$ is finite, $\K^B$ is compact,  thus $Im \varphi_{B}$ is closed, proving that it contains $\overline {\mathcal{C}_{A}}$.
 \end{proof}
 
 {\bf Claim 5} The above inclusion is an equality if in addition $\K$ is 
projectively trivial.
 
 \begin{proof} Let $g\in Im \varphi_{B}$ 
and $I'$ be a finite subset of $I$. We claim that there is some $j\in 
J$ such that $A_{-,j}$ and $g$ coincide on $I'$.  Say that two indices $j',j''\in J$ are 
equivalent if the restriction to $I'$ of the columns 
$A_{-,j'},A_{-,j''}$ are identical.  Since $I'$ and $\K$ are finite, 
this equivalence relation has only finitely many classes. Let $J'$ be 
a finite subset of $J$ containing an element of each class. The 
projection map $p:\K^{J}\rightarrow \K^{J'}$ being $1-1$ on 
$A_{I',-}:= \{A_{i,-} : i\in 
I'\}$, it is $1-1$  on $B'$, the subalgebra of $K^J$ generated by $A_{I',-}$.
Let $h\in Hom(B, \K)$ 
such that $\varphi(h)=g$. We have $g(i)=h(A_{i,-})$ for every $i\in 
I$. Since $B'$ is a subalgebra of $B$, $h$ induces an homomorphism 
from $B'$ into $\K$ and,  since $p$ is $1-1$ on $B'$, it induces an 
homomorphism $h'$ from $B''$, the image of $B'$ under $p$,  into $\K$. 
Since $\K$ is projectively trivial, $h'$ is a projection, and since, from our 
construction,  $h'(p(\overline 
{x}))=h(\overline{x})$ for every $\overline {x}\in B'$, there is some $j\in J$  such 
that $h(A_{i,-})= A(i,j)$ for all $i\in I'$. It follows that  
$g(i)=A(i,j)$ as required. 
 \end{proof}

(b) Let $f: L \rightarrow \K$ be a homomorphism from a subalgebra $L$ 
of a finite power $\K^n$ where  $\K:=\{0,1\}$ is the $2$-element Boolean 
algebra or the $2$-element lattice as well. In both cases, $f^{-1}(0)$ and 
$f^{-1}(1)$ have respectively a largest element $a:=(a_{1},\dots, 
a_{n})$  and a least element 
$b:=(b_{1},\dots, b_{n})$.  Since $b\not \leq a$ there is an indice $i$, $i<n$, such that 
$a_{i}=0$ and $b_{i}=1$. Clearly, $f$ is the restriction to $L$ of the 
$i$-th projection $p_{i}$. Hence $\K$ is projectively trivial.   
\end{proof}

\subsection{From projectively trivial algebras to algebras with the projection property}
Projectively trivial algebras seem to be interesting objects to consider, particularly in view of 
general studies about duality (as developped in \cite{clarke}). These algebras fall into a 
general class of structures, which has attracted some attention recently, those with the {\it  projection
property}. 
Let 
$\K$ be  a structure (e.g.  an algebra or a relational structure). Let $n$ be a  non-negative integer,
$\K$ has the
$n$-{\it projection property} if every idempotent homomorphism from
$\K^{n}$, the $n$-th power of $\K$,  into $\K$ is a projection. If this  holds for every $n$, we say
that $\K$ has the {\it projection property}. 

This notion was introduced  for posets by E.Corominas \cite{cor} in $1990$. Several papers have
followed (e.g. see \cite{pou-ros-sto} for reflexive relational structures, \cite {Del1}\cite{delhomme99} for reflexive graphs and \cite{larose2001} for irreflexive graphs).

Between the projectively trivial algebras and algebras with the projection property, are   algebras $\K$  which have no other  homomorphisms  from their finite powers into $\K$  than the projections. For the ease of the discussion, we will name these algebras "{\it projective}" (despite the fact that this word has an other meaning in algebra).

Clearly, an algebra $\K$ is projective if and only if it has the projection property and it  is {\it rigid} in the sense that there is no other endomorphism from $\K$ into itself other than  the identity.  A further discussion  about these algebras belongs to 
 the theory of clones.
 
 Let us recall that a  {\it clone} on $E$ is a composition closed subset of $\mathcal O$
containing
$\mathcal P$.  Equivalently, a clone is the set of term operations of some  algebra on $E$.  For example, the clone associated with  the Boolean algebra on the $2$-element set $E:=\{0,1\}$ is  $\mathcal O$, whereas the clone associated with  $2$-element lattice is the set  of all monotone operations. The structural properties of an algebra and its finite powers, namely the homomorphisms and subalgebras  are entirely determined by the clone of its term operations. 
 This is readily seen in terms of the Galois correspondence   between operations and relations defined by the  relation "$f$ preserves $\rho$". The {\it polymorph}  of a set $\mathcal G$ of relations is  the set $Pol (\mathcal G)$ of operations  which preserve  each $\rho$ in  $\mathcal G$. Polymorphs can be characterized as locally closed clones.  If $E$ is finite, there are just clones.  The {\it invariant}  of a set $\mathcal F$ of operations  is  the set $Inv(\mathcal{F})$ of relations   which are preserved by   each $f$ in  $\mathcal F$.  This is   the set of subalgebras of finite powers of $(E, \mathcal F)$.  The collection of polymorphs, once ordered by inclusion, form a complete lattice, namely  the Galois lattice of the above correspondence. If $E$ is finite, it coincides with  the lattice of clones. 
 The relation "${f}$ commutes with $g$" also defines a Galois correspondence;
  The   {\it centralizer}  of a set $\mathcal F$ of operations is the set $Z(\mathcal F)$ of operations which commute with every  $f$ in $\mathcal F$. The corresponding Galois lattice is a subset of the lattice of clones. On a finite set $E$, this lattice is finite \cite{wil-bur87}, contrarily to   the lattice of clones.

We will note the following well known fact:
  \begin{fact}\label{fact}{\it An operation $f$ is a projection (resp. is idempotent) if and only if  it commutes with all operations ( resp all constant unary operations). That is  $Z(\mathcal O)= \mathcal P$ and $ Z(\mathcal C)= \mathcal I$}.
\end{fact}

With  the notion of centralizer,  we immediately have 
\begin{fact}
 An algebra $\K:=(E, \mathcal F)$ is projective if and only if  $Z(\mathcal F)=\mathcal P$. 
 \end{fact}
 Hence, the classification of these algebras amounts to the classification   of clones whose the centralizer consists only of projections. 
Trivially, this collection of clones is a final segment of the lattice of clones  which, from Fact \ref{fact}, is non-empty. Hence, beyond $\mathcal O$, maximal clones are natural candidates. On a finite universe, there are only finitely many maximal clones, and  they have been entirely determined  \cite{rosenberg65}. A search is then possible. We do hope to report on it in the near future.

\section*{Acknowledgements} We are pleased to thank I.G.Rosenberg for his careful examination of a
preliminary version of this paper.

\end{document}